\documentclass[a4paper,11pt]{amsart}

\usepackage{amsmath}
\usepackage{amssymb}
\usepackage{amstext}
\usepackage{amsbsy}
\usepackage{amsopn}
\usepackage{amsthm}
\usepackage{amscd}
\usepackage{amsxtra}
\usepackage{amsfonts}

\usepackage{ifthen}
\usepackage{xspace}
\usepackage{fancyheadings}
\usepackage{layout}
\usepackage{hyperref}

\newcommand{\preprintserver}[2]{\href{http://xxx.lanl.gov/abs/math/#2}{#1/#2}}

\input xy
\xyoption{matrix}
\xyoption{arrow}
\xyoption{curve}
\newdir{ (}{{}*!/-5pt/@^{(}}
\newcommand{\xycenter}[1]{\begin{center}
                          \mbox{\xymatrix{#1}}
                          \end{center}
                         }

\newboolean{xlabels}
\newcommand{\xlabel}[1]{
                        \label{#1}
                        \ifthenelse{\boolean{xlabels}}
                                   {\marginpar{#1}}
                                   {}
                       }

\newcommand{\QZ}{\mathbb{Q}}
\newcommand{\ZZ}{\mathbb{Z}}

\newcommand{\PZ}{\mathbb{P}}
\newcommand{\pz}{\mathbb{P}}

\newcommand{\FZ}{\mathbb{F}}

\newcommand{\Ptwo}{{\mathbb{P}^2}}

\newcommand{\CC}{\mathbb{C}}

\newcommand{\PP}{\mathbb{P}}

\newcommand{\FF}{\mathbb{F}}

\newcommand{\sF}{{\mathcal F}}
\newcommand{\sG}{{\mathcal G}}

\newcommand{\sO}{{\mathcal O}}

\newcommand{\suchthat}{\, | \,}

\newboolean{probleme}
\newcommand{\problem}[1]
           {\ifthenelse{\boolean{probleme}}
                       {{\bf(PROBLEM: #1)\bf}}
                       {}
           }
\newboolean{zukuenftiges}
\newcommand{\zukunft}[1]
           {\ifthenelse{\boolean{zukuenftiges}}
                       {{\bf(AUSBAUM\"OGLICHKEIT: #1)\bf}}
                       {}
           }
\newboolean{extras}
\newcommand{\extra}[1]
           {\ifthenelse{\boolean{extras}}
                       {{\bf EXTRA #1 EXTRA\bf}}
                       {}
           }

\newboolean{ignore}
\newcommand{\ignore}[1]
           {\ifthenelse{\boolean{ignore}}
                       {{\bf IGNORE #1 IGNORE\bf}}
                       {}
           }

\DeclareMathOperator{\Hilb}{Hilb}

\DeclareMathOperator{\spec}{spec}

\theoremstyle{plain}
\begingroup

\newtheorem{thm}{Theorem}
\newtheorem{cor}[thm]{Corollary}
\newtheorem{lem}[thm]{Lemma}

\newtheorem{prop}[thm]{Proposition}

\numberwithin{thm}{subsection} 
\endgroup

\newtheorem*{thm*}{Theorem}
\newtheorem*{conj*}{Conjecture}
\newtheorem*{verm*}{Vermutung}

\theoremstyle{definition}

\newtheorem{rem}[thm]{Remark}
\newtheorem{example}[thm]{Example}

\newtheorem{script}[thm]{Script}
\numberwithin{equation}{section}

\newcommand{\nosubsections}{\renewcommand{\thethm}{\thesection.\arabic{thm}}
                            \setcounter{thm}{0}
                           }

\newcommand{\cref}[3]{(\ref{#1}, #2 \ref{#3})}

\date{\today}

\usepackage{amssymb,amsmath}
\usepackage{verbatim}
\usepackage{graphicx}
\setboolean{probleme}{true}

\newcommand{\secemail}{
\setlength{\unitlength}{1pt}
bothmer
\begin{picture}(0,1)
\put(0,0){m}
\put(-5,0){@}
\end{picture}
ath.uni-hannover.de}

\begin{document}

\title{A new family of rational surfaces in $\PZ^4$}

\address{Institiut f\"ur Mathematik\\
          Universit\"at Hannover\\
          Welfengarten 1\\
          D-30167 Hannnover
         }

\email{\secemail}

\urladdr{\href{http://www-ifm.math.uni-hannover.de/~bothmer}{http://www-ifm.math.uni-hannover.de/\textasciitilde bothmer}}

\thanks{Supported by the Schwerpunktprogramm ``Global Methods in Complex
        Geometry'' of the Deutsche Forschungs Gemeinschaft.}

\author[v. Bothmer]{Hans-Christian Graf v. Bothmer}
\author[Erdenberger]{Cord Erdenberger}
\author[Ludwig]{Katharina Ludwig}

\begin{abstract}
We describe a new method of constructing rational surfaces with given invariants in $\PP^4$ 
and present a family of degree $11$ rational surfaces of sectional genus $11$ with 
$2$ six-secants that we found with this method.
\end{abstract}

\maketitle

\section{Introduction}
\nosubsections

In 1989 Ellingsrud and Peskine showed that the degree of non-general type surfaces in $\PZ^4$
is bounded \cite{EllingsrudPeskine}. Since then their degree bound has been sharpened by various
authors, most recently by Decker and Schreyer \cite{DeckerSchreyer} to $52$. On the other hand
one has tried to construct and classify non-general type surfaces in $\PZ^4$. \cite{DES} lists the $51$ 
families of such surfaces known at that time of which $18$ are rational. Since then \cite{smallFields} 
has found $4$ more families, $3$ of them parameterizing rational surfaces. Five more families of
non-rational, non-general type surfaces were found by \cite{Aure}, \cite{AboDegree8} and 
\cite{AboDegree12}. Recently Abo announced the existence of a family of degree $12$ rational surfaces 
in $\PZ^4$. Non-general type surfaces are classified up to degree $10$ (see \cite{DeckerSchreyer} for 
an overview and references), the largest known degree of a non-general type surface in $\PZ^4$ is $15$. 
Rational surfaces are only known up to degree $12$.

In this paper we describe a new method for finding rational surfaces in $\PZ^4$ and present a new 
family of rational degree $11$ surfaces in $\PZ^4$ which we found with this method.

Our method is partly based on an idea of Schreyer \cite{smallFields} who explicitly constructed 
surfaces in $\PZ^4$ over small fields using computer algebra programs and provided a method of lifting 
these
surfaces to characteristic $0$. There he used the observation that modules with special syzygies are
much more common over small fields than over characteristic $0$ to find those modules by a random 
search. Here we use a random search over $\FZ_2$ to find linear systems with special configurations of 
basepoints on $\PZ^2$.

We construct one of our new surfaces over $\FZ_2$ in Section \ref{surface}. In Section \ref{lift} we 
prove that this surface lies in a family that is also defined in characteristic zero. Finally we  
explain our search-algorithm in Section \ref{search}.

Part of our calculations were done at the Gauss Laboratory at the University of G\"ottingen. We would like to thank Yuri Tschinkel for giving us this opportunity.
\section{The surface} \label{surface}
\nosubsections

\newcommand{\FFtwo}{{\FZ_2}}
\renewcommand{\FF}[1]{\FZ_{2^{#1}}}

Let us first fix some notation. We work over the fields $\FFtwo$, $\FF{14}$ and $\FF{5}$ which we
realize as
\[
    \FF{14} = \FFtwo[t]/(t^{14}+t^{13}+t^{11}+t^{10}+t^{8}+t^{6}+t^{4}+{t}+1)
\]
and $\FF{5} = \FFtwo[t]/(t^{5}+t^{3}+t^{2}+{t}+1)$.

Over these fields we consider the points
\begin{align*}
    P &= (0 : 0 : 1) \in \PZ^2(\FFtwo)\\
    Q &= (t^{11898}, t^{137}, 1) \in \PZ^2(\FF{14})\\
    R &= (t^6 : t^{15} : 1) \in \PZ^2(\FF{5}).
\end{align*}

\begin{lem} \label{l-orbits}
The orbits of $Q$ and $R$ under the Frobenius-endomorphism are of degree $14$ and $5$
respectively. We denote the corresponding points by $Q_1,\dots,Q_{14}$ and $R_1,\dots,R_5$.
\end{lem}

\begin{proof}
The orbits of $Q$ and $R$ are defined by the kernels of
\[
	\FFtwo[x,y,z] \to \FF{14}[x,y,z]/I_Q
\]
and
\[
	\FFtwo[x,y,z] \to \FF{5}[x,y,z]/I_R
\]
where $I_Q$ and $I_R$ are the ideals of $Q$ and $R$ respectively. Using
Script \ref{script} one can calculate these kernels and check that the degrees of their vanishing sets
are $14$ and $5$.
\end{proof}

\begin{prop} \label{p-rank5}
Let $L$ be the class of a line in $\PZ^2$. Then
\[
       |9L-3P-2Q_1-\dots-2Q_{14}-R_1-\dots-R_5| = \PZ^4
\]
and this linear system has only the assigned base points.
\end{prop}

\begin{proof}
Using Script \ref{script} we intersect the ideal of $3P$, the ideal of the orbit of $2Q$ and
the ideal of the orbit of $R$. We check that the intersection contains exactly $5$ independent $9$-tics whose baselocus is of degree $1\cdot 6 + 14 \cdot 3 + 5 \cdot 1 = 53$.
\end{proof}

\begin{prop} \label{p-smooth}
Let $S$ be the blowup of $\PZ^2$ in $P, Q_1,\dots,Q_{14}, R_1, \dots, R_{5}$ and $E_1,\dots,E_{20}$ be 
the corresponding exceptional divisors. Then the
linear system
\[
    \left| 9L - 3E_1 - \sum_{i=2}^{15} 2E_i - \sum_{i=16}^{20} E_i\right|
\]
is very ample and embeds $S \subset \PZ^4$ as a smooth surface of degree $11$ and sectional
genus $11$.
\end{prop}

\begin{proof}
In Script \ref{script} we define a morphism
\[
	\FFtwo[x_0,\dots,x_4] \to \FFtwo[x,y,z]
\]
using the $5$ independent $9$-tics found. The kernel of this map is the ideal of $S \subset \PP^4$.
We then calculate that $S$ has degree $11$, and sectional genus $11$ and finally check smoothness by the Jacobi criterion.
\end{proof}

\begin{prop} \label{p-secants}
$S \subset \PZ^4$ has two $6$-secants. In particular $S$ can not lie in one of the
known families.
\end{prop}

\begin{proof}
Every $6$-secant of $S$ must be contained in all quintics that contain $S$. The final lines of Script \ref{script}
calculate that the vanishing locus $(I_S)_{\le 5}$ contains $S$ and two lines, which turn out to
be $6$-secants. The known families of rational surfaces of degree $11$ and sectional genus $11$
have $0$, $1$ or infinitely many $6$-secants \cite{DES}. 
\end{proof}

\section{Lifting to characteristic zero} \label{lift}
\nosubsections

\newcommand{\Fp}{{\FZ_p}}

In this section we denote schemes defined over $\spec \ZZ$ with a subscript $\ZZ$ and their fibers over
points of $\spec \ZZ$ with subscripts $\Fp$ or $\QZ$.

\begin{prop}[Schreyer \cite{smallFields}] \label{lifting}
Consider a smooth projective variety $X_\ZZ \subset \PZ^N_\ZZ$ over $\spec \ZZ$, a map
\[
    \phi \colon \sF \to \sG
\]
between vector bundles of rank $f$ and $g$ on $X_\ZZ$, the determinantal
subvariety $Y_\ZZ \subset X_\ZZ$ where $\phi$ has rank $k$, and a $\FZ_p$-rational point $y \in Y_\Fp$.

If the tangent space $T_{Y_\Fp,y}$
of $Y_\Fp$ in $y$ is a linear subspace of codimension $(f-k)(g-k)$ in $T_{X_\Fp,y}$ then $y$ lies on an 
irreducible component $Z_\ZZ$ of $Y_\ZZ$ that has nonempty fibers over an open subscheme of $\spec 
\ZZ$.
\end{prop}

\begin{proof}
Since $Y_\Fp$ is determinantal, the codimension of $Y_\Fp$ in $X_\Fp$ is at most $c = (f-k)(g-k)$. The 
condition on the tangent space ensures that $Y_\Fp$ is smooth of this codimension in $y$, or 
equivalently $Y_\Fp$ is of dimension $d=\dim X_\Fp-c$ in $y$. Let now $Z_\ZZ$ be a component of $Y_\ZZ$ 
that contains $y$. Since $Y_\ZZ$ is determinantal in $X_\ZZ$ each component of $Y_\ZZ$ has also at most 
codimension $c$ \cite[Ex 10.9, p. 246]{Ei95}.
Since $\dim X_\ZZ = \dim X_\Fp + 1$ the dimension of $Z_\ZZ$ is at
least $d+1$. $Z_\Fp$ contains $y$ and is therefore of dimension at most $d$. Hence $Z_\ZZ$ can
not be contained in the fiber $Y_\Fp$ and has nonempty fibers over an open subscheme of $\spec \ZZ$.
\end{proof}

\newcommand{\PtwoZZ}{{\PZ^2_\ZZ}}
\newcommand{\sOPtwoZZ}{{\sO_\PtwoZZ}}

On $\PtwoZZ$ we have the map
\[
    \tau_k \colon H^0(\sOPtwoZZ(a)) \to \sOPtwoZZ(a) \oplus
                                3 \sOPtwoZZ(a-1) \oplus
                                \dots \oplus
                                { k+2 \choose 2} \sOPtwoZZ(a-k)
\]
that associates to each polynomial of degree $a$ the coefficients of its taylor expansion up to degree 
$k$.

\begin{lem}
If $a > k$ then the image of $\tau_k$ is a vector bundle $\sF_k$ of rank ${k+2 \choose 2}$ over $\spec 
\ZZ$.
\end{lem}

\begin{proof}
In each point we consider an affine $2$-dimensional neighborhood where we can choose the
${k+2 \choose 2}$ coefficients of the affine taylor expansion independently. This shows that
the image has at least this rank everywhere. That this is also the maximal rank follows from the Euler 
relation.
\end{proof}

\begin{rem}
Notice that the morphism $H^0(\sOPtwoZZ(a)) \to { k+2 \choose 2} \sOPtwoZZ(a-k)$ is not
surjective in characteristics that divide $a$. One really has to consider the whole taylor expansion.
\end{rem}

Set now $X_\ZZ = \Hilb_{1,\ZZ} \times \Hilb_{14,\ZZ} \times \Hilb_{5,\ZZ}$ where $\Hilb_{k,\ZZ}$ 
denotes the Hilbert scheme of $k$ points in $\PP^2_\ZZ$ over $\spec \ZZ$, and let
\[
    Y_\ZZ = \{ (p,q,r) \suchthat h^0(9L - 3p - 2q - 1r) \ge 5\} \subset X_\ZZ
\]
be the subset where the linear system of ninetics with triplepoint in $p$, doublepoints in $q$ and
single basepoints in $r$ is at least of projective dimension $4$.

\begin{prop}
There exist vector bundles $\sF$ and $\sG$ of ranks $55$ and $53$ respectively on $X_\ZZ$ and a 
morphism
\[
    \phi  \colon \sF \otimes \sO_{X_\ZZ} \to \sG
\]
such that the determinantal locus where $\phi$ has rank $50$ is supported on $Y_\ZZ$.
\end{prop}

\begin{proof}
On the cartesian product
\xycenter{
     \Hilb_{d,\ZZ} \times \PtwoZZ \ar[r]^-{\pi_2} \ar[d]^-{\pi_1} & \PtwoZZ \\
      \Hilb_{d,\ZZ}
            }
we have the morphisms
\[
    \pi_2^*\tau_k \colon H^0(\sOPtwoZZ(9))\otimes \sO_{\Hilb_{d,\ZZ} \times \PtwoZZ} \to \pi_2^* \sF_k.
\]
Let now $P_d \subset  \Hilb_{d,\ZZ} \times \PtwoZZ$ be the universal set of points. Then
$P_d$ is a flat family of degree $d$ over $\Hilb_{d,\ZZ}$ and
\[
    \sG_{k} := (\pi_1)_* ((\pi_2^*\sF_k)|_{P_d})
\]
is a vector bundle of rank $d{ k+2 \choose 2}$ over $\Hilb_{d,\ZZ}$. On
\[
      X_\ZZ = \Hilb_{1,\ZZ} \times \Hilb_{14,\ZZ} \times \Hilb_{5,\ZZ}
\]
the induced map
\[
    \phi \colon H^0(\sOPtwoZZ(9))\otimes \sO_{X_\ZZ} \xrightarrow{\tau_2\oplus\tau_1\oplus\tau_0}
                         \sigma_1^* \sG_{2} \oplus \sigma_{14} ^* \sG_{1}
                         \oplus \sigma_{5}^* \sG_{0}
\]
has the desired properties, where $\sigma_d$ denotes the projection to $\Hilb_{d,\ZZ}$.
\end{proof}

\begin{thm}
There exists a family of smooth rational surfaces in $\PZ^4_\CC$ with $d=11$, $\pi = 11$, $K^2 = -11$
and two $6$-secants.
\end{thm}

\begin{proof}
By determining the infinitesimal deformations of our $20$ points $P,Q_1,\dots,Q_{14},R_1, \dots ,R_5$
in $\PZ^2_\FFtwo$ we can check with a Macaulay calculation that
\[
    y := (P,\{Q_1,\dots,Q_{14}\},\{R_1,\dots,R_5\}) \in Y_\FFtwo \subset X_\FFtwo
\]
satisfies the conditions of Proposition \ref{lifting}. A script performing these calculations can be obtained 
from our webpage \cite{ratsurfweb}. We therefore have a component $Z_\ZZ$ of $Y_\ZZ$
that contains our configuration of basepoints. Since the conditions
\begin{enumerate}
\item the points of $p$, $q$ and $r$ are distinct
\item $h^0(9L-3p-2q-1r)=5$
\item the linear system $|9L-3p-2q-1r|$ has no further basepoints
\item the image of the corresponding rational map $\phi \colon \Ptwo \to \PZ^4$ is a smooth surface $S$
\item $S$ has two $6$-secants
\end{enumerate}
are all open on $Z_\ZZ$ and $y$ is a point on this component that satisfies all conditions, they must 
hold on a nonempty open subset of $Z_\ZZ$. Since  $Z_\ZZ$ is irreducible and has nonzero fibers over 
the generic point, we obtain smooth surfaces in characteristic zero. The invariants can be calculated 
from the multiplicities of the baselocus using the following proposition.
\end{proof}

\begin{prop} \label{p-formulae1}
Let $S = \PZ^2_\CC(p_1,\dots,p_l)$ be the blowup of $\PZ^2_\CC$ in $l$ distinct points. We denote by 
$E_1,\dots,E_l$ the corresponding exceptional divisors and by $L$ the pullback of a general line in 
$\PP^2_\CC$ to $S$. If
$
    | aL - \sum_{i=1}^l b_i E_i |
$
is a very ample linear system of dimension $4$ for suitable $a$ and $b_i$, then $S \subset \PZ^4_\CC$ 
is a rational surface of degree
\[
    d = a^2 - \sum_{i=1}^l b_i^2
\]
and sectional genus
\[
    \pi = { a-1 \choose 2} - \sum_{i=1}^l { b_i \choose 2}.
\]
The self-intersection of the canonical divisor of $S$ is $K^2 = 9-l$.
\end{prop}

\begin{proof}
Set $H = aL - \sum_{i=1}^l b_i E_i$. Then
\[
    d = H^2 = (aL - \sum_{i=1}^l b_i E_i)^2 = a^2 - \sum_{i=1}^l b_i^2
\]
since $L^2=1$, $L\cdot E_i=0$ and $E_i \cdot E_j = -\delta_{ij}$.
The canonical divisor of $S$ is $K = -3L+\sum_{i=1}^l E_i$, so $K^2=9-l$. The sectional genus
of $S$ can be calculated by adjunction:
\begin{align*}
    \pi &= \frac{1}{2}H(K+H) + 1 
         = { a-1 \choose 2} - \sum_{i=1}^l { b_i \choose 2}.
\end{align*}
\end{proof}

\section{The search} \label{search}
\nosubsections

In this section we will describe our search-algorithm. We first need to find suitable linear systems 
for given invariants. For that we make the following observation which is a direct consequence of 
Proposition \ref{p-formulae1}:

\begin{cor}
In the situation of Proposition \ref{p-formulae1} we set $\beta_j = \#\{i \suchthat b_i=j\}$. The 
invariants
of $S$ are then linear forms in the $\beta_j$'s:
\begin{align*}
    d     &= a^2 - \sum_j \beta_j j^2\\
    \pi   &=  { a-1 \choose 2} - \sum_{j}\beta_j { j\choose 2}\\
    K^2 &= 9 - \sum_j \beta_j.
\end{align*}
\end{cor}

For given $d$, $\pi$, $K^2$ and $a$  the linear system above has only finitely many integer solutions.
One can find these solutions by integer programming. We have used an algorithm from \cite[Chapter 
8]{CLOusingAG}.

\begin{example}
For $d=11$, $\pi=11$, $K^2=-11$ and $a=9$ the only solution is $\beta_3=1$, $\beta_2=14$ and 
$\beta_1=5$.
\end{example}

For a given set of $\beta_j$'s we have chosen random points in $\PZ^2$ over $\FF{\beta_j}$, checked
if their orbit under the Frobenius endomorphism had degree $\beta_j$; checked whether the corresponding 
linear system
$
    | aL - \sum_{i=1}^l b_i E_i |
$
was $4$-dimensional; calculated the image of the corresponding map to $\PZ^4$ and checked whether
this image was a smooth rational surface. The script we used is available on our web page \cite{ratsurfweb}.

\begin{figure}[b]
\begin{center}
\begin{tabular}{|l|r|r|r|r|r|r|r|}
\hline
Example & $d$ & $\pi$ & speciality & trials & surfaces & rate & log rate \\ \hline
B1.7 & 5 & 2 & 0 & 1000 & 871 & 87,1\% & -0,2\\ \hline
B1.8 & 6 & 3 & 0 & 1000 & 311 & 31,1\% & -1,7\\ \hline
B1.9 & 7 & 4 & 0 & 1000 & 188 & 18,8\% & -2,4\\ \hline
B1.10 & 8 & 5 & 0 & 1000 & 312 & 31,2\% & -1,7\\ \hline
B1.11 & 8 & 6 & 1 & 10000 & 184 & 1,84\% & -5,8\\ \hline
B1.12 & 9 & 6 & 0 & 10000 & 2173 & 21,73\% & -2,2\\ \hline
B1.13 & 9 & 7 & 1 & 100000 & 446 & 0,446\% & -7,8\\ \hline
B1.14 & 10 & 8 & 1 & 1000000 & 0 & 0,0000\% & $-\infty$\\ \hline
B1.15 & 10 & 9 & 2 & 1000000 & 42 & 0,0042\% & -14,5\\ \hline
B1.16 & 10 & 9 & 2 & 1000000 & 267 & 0,0267\% & -11,9\\ \hline
B1.17 & 11 & 11 & 3 & 10000000 & 0 & 0,00000\% & $-\infty$\\ \hline
B1.18 & 11 & 11 & 3 & 10000000 & 5 & 0,00005\% & -20,9\\ \hline
B1.19 & 11 & 11 & 3 & 10000000 & 0 & 0,00000\% & $-\infty$\\ \hline
New & 11 & 11 & 3 & 2000000 & 21 & 0,00105\% & -16,5\\ \hline
\end{tabular}
\caption{Results of random searches using our script. The numbering of the examples
is as in \cite{DES}.} \xlabel{f-results}
\end{center}
\end{figure}

\begin{figure}
\includegraphics[width=12.5cm, trim=0 0 0 0cm]
                              {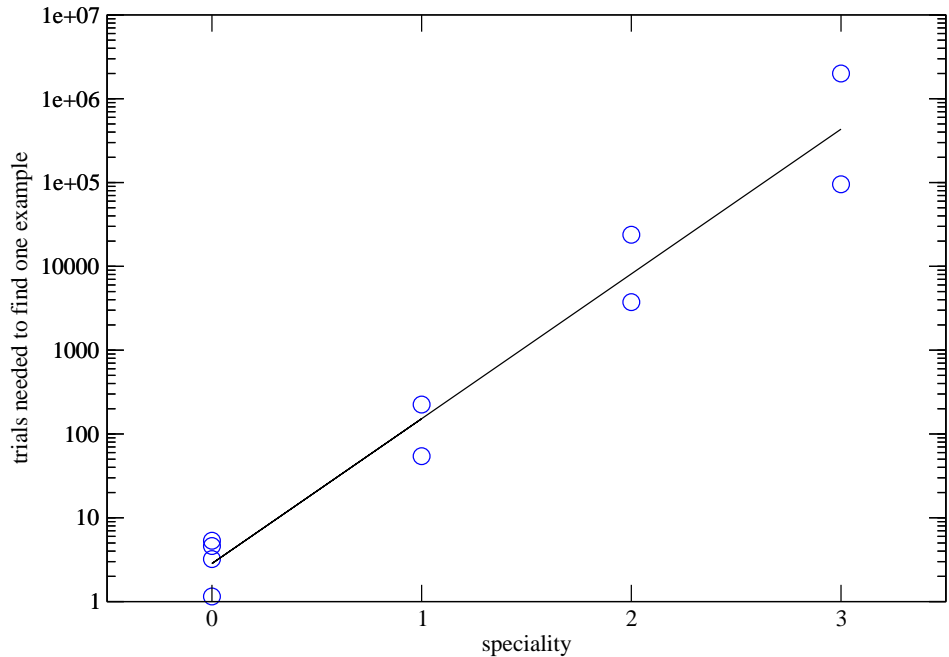}
\caption{The difficulty of finding a surface grows exponentially with the speciality.}
\xlabel{f-trials}                              
\end{figure}

For comparison we also tried to reconstruct the rational surfaces of \cite{DES} using the
basepoint multiplicities provided there. Our results are collected in Figure \ref{f-results} and Figure \ref{f-trials}. Notice that the number of trials needed to find an example grows exponentially with the speciality as expected. From this we expect that surfaces of speciality $4$ can be found in approximately $500$ Million trials. Using our program this would take about $500$ weeks on a $2$ GHz machine.

\begin{rem}
Notice that our approach can only find linear systems where the points in each group of constant 
multiplicity are in uniform position.
\end{rem}

\begin{rem}
Another way of constructing random groups of points is via syzygies and the Theorem of Hilbert-Birch.
We have also tried this, but found the above approach more effective.
\end{rem}

\section{Appendix}
\nosubsections

Here we provide a script for the computer algebra program Macaulay 2 \cite{M2} that does the calculations needed in Section \ref{surface}. The script
can also be obtained from our webpage \cite{ratsurfweb}.

\begin{script} \quad \xlabel{script}

\begin{small}
\begin{verbatim}
-- construct a surface over F_2 using frobenius orbits

-- define coordinate ring of P^2 over F_2
F2 = GF(2)
S2  = F2[x,y,z]

-- define coordinate ring of P^2 over F_2^14 and F_2^5
St  = F2[x,y,z,t]
use St; I14 = ideal(t^14+t^13+t^11+t^10+t^8+t^6+t^4+t+1); S14 = St/I14
use St; I5 = ideal(t^5+t^3+t^2+t+1); S5 = St/I5

-- the points
use S2; P = matrix{{0_S2, 0_S2, 1_S2}}
use S14;Q = matrix{{t^11898, t^137, 1_S14}}
use S5; R = matrix{{t^6, t^15, 1_S5}}

-- their ideals
IP = ideal ((vars S2)*syz P)
IQ = ideal ((vars S14)_{0..2}*syz Q)
IR = ideal ((vars S5)_{0..2}*syz R)

-- their orbits
f14 = map(S14/IQ,S2); Qorbit = ker f14
degree Qorbit   -- degree = 14

f5 = map(S5/IR,S2); Rorbit = ker f5
degree Rorbit   -- degree = 5

-- ideal of 3P
P3 = IP^3;

-- orbit of 2Q
f14square = map(S14/IQ^2,S2); Q2orbit = ker f14square;

-- ideal of 3P + 2Qorbit + 1Rorbit
I = intersect(P3,Q2orbit,Rorbit);

-- extract 9-tics
H = super basis(9,I)
rank source H   -- affine dimension = 5

-- count basepoints (with multiplicities)
degree ideal H   -- degree = 53

-- construct map to P^4
T = F2[x0,x1,x2,x3,x4]
fH = map(S2,T,H);

-- calculate the ideal of the image
Isurface = ker fH; 

-- check invariants
betti res coker gens Isurface
codim Isurface    -- codim = 2
degree Isurface   -- degree = 11
genera Isurface   -- genera = {0,11,10}

-- check smoothness
J = jacobian Isurface;
mJ = minors(2,J) + Isurface;
codim mJ  -- codim = 5

-- count 6-secants
-- ideal of 1 quartic and 5 quintics
Iquintics = ideal (mingens Isurface)_{0..5}; 

-- calculate the extra components where these vanish
secants = Iquintics : Isurface;
codim secants  -- codim = 3
degree secants  -- degree = 2
secantlist = decompose secants  -- two components

-- check number of intersections
degree (Isurface+secantlist#0)  -- degree = 6
codim (Isurface+secantlist#0)  -- codim = 4
degree (Isurface+secantlist#1)  -- degree = 6
codim (Isurface+secantlist#1)  -- codim = 4
\end{verbatim}
\end{small}
\end{script}



\newcommand{\etalchar}[1]{$^{#1}$}
\def\cprime{$'$}

\end{document}